\documentclass[12pt, a4paper]{article}
\usepackage[pagewise]{lineno}
\usepackage{graphicx}
\usepackage{amssymb}
\usepackage{latexsym, bm}
\usepackage{multicol}
\usepackage{indentfirst}
\usepackage{amssymb,amsfonts}
\usepackage{amsmath}
\usepackage{cases}
\usepackage[noadjust]{cite}
\usepackage[colorlinks,citecolor=red]{hyperref}
\newtheorem{theorem}{Theorem}
\newtheorem{lemma}{Lemma}

\newtheorem{proposition}{Proposition}
\newtheorem{definition}{Definition}

\usepackage{amssymb}
\numberwithin{equation}{section}
\numberwithin{theorem}{section}
\numberwithin{remark}{section}
\numberwithin{definition}{section}
\numberwithin{lemma}{section}
\numberwithin{corollary}{section}
\numberwithin{proposition}{section}
\numberwithin{notation}{section}
\title{Liouville-type theorem for the stationary inhomogeneous Navier-Stokes equations}
\author{Huiting Ding\footnote{dinghuiting27@163.com}\quad Wenke Tan\footnote{tanwenkeybfq@163.com}\\
{\small Key Laboratory of Computing and Stochastic Mathematics (Ministry of Education),}\\
{\small School of Mathematics and Statistics, Hunan Normal University,}\\
{\small Changsha, Hunan 410081, China}\\
}

\date{}

\begin{document}
\maketitle
{\bf Abstract:} In this manuscript, a new Liouville-type theorem for the three-dimensional stationary inhomogeneous Navier-Stokes equations is established. We first
localize the Dirichlet energy into the region near the origin in frequency spaces by two times localizations. The first localization is to eliminate the non-zero frequency part coming from the interaction between $\rho u$ and $u$, the second one is to eliminate the non-zero frequency part coming from the interaction between $\rho$ and $u$. Based on the local formula of Dirichlet energy, we can establish suitable estimates on different frequency parts of $u$ and $\rho$, then show our new Liouville-type theorem.

\medskip
{\bf Mathematics Subject Classification (2020):} \  35Q30, 76B03, 35Q35.
\medskip

{\bf Keywords:}  Stationary inhomogeneous Navier-Stokes equations; Liouville type theorem; Frequency spaces
\section{Introduction}
In this paper, we consider the following stationary incompressible inhomogeneous Navier-Stokes system in the whole space $\mathbb{R}^3$:
\begin{equation}\label{NS1}
 \left\{\begin{array}{ll}
div(\rho u)=0,\\
div(\rho u\otimes u)-\mu\nabla u+\nabla P=0,\\
\nabla\cdot u=0.
\end{array}\right.
\end{equation}
Here $\rho=\rho(x)\in \mathbb{R}^+$ and $u=u(x)\in \mathbb{R}^3$ stand for the density and velocity field respectively, and $P$ is a scalar pressure function. The viscosity coefficient $\mu>0$ is a given positive real number. Throughout, we assume $\mu=1$ for convenience.
\par
From the mathematical point of view, the stationary incompressible inhomogeneous Navier-Stokes system \eqref{NS1} is significantly related to following stationary incompressible homogeneous Navier-Stokes equations in $\mathbb{R}^3$:
\begin{equation}\label{NS}
 \left\{\begin{array}{ll}
-\Delta u+u\cdot\nabla u+\nabla P=0,\\
\nabla\cdot u=0.\\
\end{array}\right.
\end{equation}
It is well known that due to the pioneering work of Leray \cite{Leray}, it has been an open problem whether $u = 0$ is the only solution with finite Dirichlet integral
\begin{align}\label{Du}
\int_{\mathbb{R}^3}|\nabla u|^2dx<+\infty.
\end{align}
This is a famous Liouville type statement on the stationary Navier-Stokes equations.
Galdi [\cite{Galdi},
Theorem X.9.5] proved that if $u\in L^{\frac{9}{2}}(\mathbb{R}^3)$, then $u=0$. The
same proof works for dimension $n\geq4$  by assuming the smooth solution $u$ obeys \eqref{Du} and additional assumption $\lim\limits_{|x|\to\infty}|u(x)|=0$.
This result was improved by a log factor in Chae and Wolf \cite{Chae} by assuming that
\begin{align*}
\int_{\mathbb{R}^3}|u|^{\frac{9}{2}}\left\{\log\left(2+\frac{1}{|u|}\right)\right\}^{-1}dx<\infty.
\end{align*}
Chae \cite{Chae2} showed that the condition $\Delta u\in L^{\frac{6}{5}}(\mathbb{R}^3)$ is sufficient for $u=0$ in $\mathbb{R}^3$.
Seregin \cite{Seregin2} proved that a solution in $\mathbb{R}^3$ is $0$ if $u\in L^6(\mathbb{R}^3)\cap BMO^{-1}(\mathbb{R}^3)$. Kozono et al. \cite{kt} showed the condition
\begin{align*}
\limsup\limits_{|x|\to\infty}|x|^{\frac{5}{3}}|\omega(x)|\leq(\delta D(u))^{\frac{1}{3}}
\end{align*}
or
\begin{align*}
\|u\|_{L^{\frac{9}{2},\infty}(\mathbb{R}^3)}\leq(\delta D(u))^{\frac{1}{3}},
\end{align*}
for a small constant $\delta>0$, is sufficient for $u=0$ in $\mathbb{R}^3$. Then, the restriction imposed on the norm $\|u\|_{L^{\frac{9}{2},\infty}(\mathbb{R}^3)}$ was relaxed by Seregin and Wang in \cite{Seregin3}.
{Chamorro, Jarr\'{\i}n and Lemari\'{e}-Rieusset \cite{Chamorro} proved $u=0$ by assuming that
$$u\in L^p(\mathbb{R}^3) \quad with \quad 3\leq p\leq\frac{9}{2}$$
or
$$u\in L^p(\mathbb{R}^3)\cap \dot{B}^{\frac{3}{p}-\frac{3}{2}}_{\infty,\infty}(\mathbb{R}^3)
\quad with\quad \frac{9}{2}<p<6.$$
Yuan and Xiao \cite{Yuan} improved the lower bound for $p$ so that they considered $2\leq p <\frac{18}{5}$. Moreover, Jarr\'{\i}n \cite{Jarrin} extended previous result in \cite{Chamorro} to the $L^p$-spaces with $\frac{3}{2}<p<3$ or $L^p\cap\dot{H}^{-1}$  with $\frac{9}{2}< p<+\infty$.
For more references on Liouville type theorem for stationary Navier-Stokes equations, we refer to \cite{Tsai,Cho,Chae3,Chae4} and the references therein.
Recently, Tan \cite{Tan} showed that the global defined quantity \eqref{Du} is completely determined by the local information of the velocity $u$ near the origin in frequency spaces and established a new Liouville-type theorem under the assumption
$\lim\limits_{k\rightarrow-\infty}\|\dot{S}_ku\|_{\dot{B}^{-1}_{\infty,\infty}(\mathbb{R}^3)}<\infty.$}

To the best of our knowledge, there are relatively few studies on the Liouville property of the steady-state
solution for the incompressible inhomogeneous Navier-Stokes equations. In \cite{Li-Yu}, Li and Yu pointed out that the condition
\begin{align}\label{li-yu}
\|\rho\|_{L^{\infty}(\mathbb{R}^3)}+\|\nabla u\|_{L^2(\mathbb{R}^3)}<\infty
\end{align}
is very natural because most physical flows have bounded density and finite enstrophy. Liu and Liu \cite{Liu-Liu} proved taht $u=0$ if
$$\|\rho\|_{L^{\infty}(\mathbb{R}^3)}+\|u\|_{L^p(\mathbb{R}^3)}<\infty \quad and\quad \lim\limits_{|x|\to\infty}|u(x)|=0$$
with $p\in[1,\frac{9}{4})$.
They also extended their result to Lorentz spaces by assuming
$$\|\rho\|_{L^{\infty}(\mathbb{R}^3)}+\|u\|_{L^{p,q}(\mathbb{R}^3)}<\infty$$
with $(p,q)\in (2,\frac{9}{2})\times[1,\infty]$ or $p = q = 2$.
Kim \cite{Kim} established the Liouville type theorem using the assumption
$$\|\rho\|_{L^{\infty}(\mathbb{R}^3)}+\|\nabla u\|_{L^2(\mathbb{R}^3)}+\|u\|_{L^{\frac{9}{2},q}(\mathbb{R}^3)}<\infty,$$
where $\frac{9}{2}\leq q<\infty$. Recently, Liu \cite{Liu} demonstrated that if
\begin{align*}
\|\rho\|_{L^{\infty}(\mathbb{R}^3)}+\|u\|_{BMO^{-1}(\mathbb{R}^3)}<\infty\quad and\quad \lim\limits_{|x|\to\infty}|u(x)|=0,
\end{align*}
then the solution is trivial.

To the authors knowledge, those authors of the previous works in the literature are analyzing the Liouville problem in physical space. There are not related results that consider the Liouville problem in the frequency space. {Inspired by the work of Tan \cite{Tan}, the purpose of this paper is to study the Liouville-type problem for the equations \eqref{NS1} in the frequency space and show that the uniqueness of weak solutions to \eqref{NS1} with assumptions \eqref{li-yu} is completely determined by the local information of $u$ near the origin in frequency spaces.}
\par
{We state our result as follows:}

\begin{theorem}\label{th1}
Assume $(\rho,u)$ be a smooth solution of \eqref{NS1} with $\rho\in L^{\infty}(\mathbb{R}^3)$ and $\nabla u\in L^2(\mathbb{R}^3)$. If $u$ satisfies the following condition
\begin{align}\label{con1}
\liminf\limits_{k\rightarrow-\infty}\|\dot{S}_{k} u\|_{\dot{B}^{-1}_{\infty,\infty}}<\infty.
\end{align}
Then it holds that $u=0$ in $\mathbb{R}^3$.
\end{theorem}



Our paper is organized as follows: In Section 2, we list some notations and recall some Lemmas about the Littlewood-Paley theory which will be used in the sequel. In Section 3, we present the proofs of results.

\section{Preliminaries}

In this section, we introduce some notations and useful auxiliary lemmas which will play a fundamental role in this paper.
\par
Firstly, let us recall that for all $f\in \mathcal{S}'$, the Fourier transform $\hat{ f}$, is defined by
\begin{align*}
\mathcal{F}(f)(\xi)=\hat{ f}(\xi)=
{(2\pi)^{-\frac{3}{2}}}\int_{\mathbb{R}^3}e^{-ix\xi}f(x)dx\ \ for\ any\ \xi\in\mathbb{R}^3.
\end{align*}
The inverse Fourier transform of any $g$ is given by
\begin{align*}
\mathcal{F}^{-1}(g)(x)=\check{g}(x)={(2\pi)^{-\frac{3}{2}}}\int_{\mathbb{R}^3}e^{ix\xi}g(\xi)d\xi\ \ for\ any\ x\in{\mathbb{R}^3.}
\end{align*}
Next, we will recall some facts about the Littlewood-Paley(L-P) decomposition of distributions into dyadic blocks of frequencies:
\par
Let $\varphi\in\mathcal{C}^{\infty}_c(\mathbb{R}^3)$ be a non negative function so that $\varphi(\xi)=1$ for $|\xi|\leq\frac{1}{2}$ and $\varphi(\xi)=0$ for $|\xi|\geq1$. Let $\psi=\varphi(\frac{\xi}{2})-\varphi(\xi)$. For every $u\in \mathcal{S}'(\mathbb{R}^3)$, the homogeneous dyadic blocks $\dot{\Delta}_j$ and homogeneous low-frequency cut-off operators $\dot{S}_j$ are defined by
\begin{align*}
&\dot{\Delta}_ju:=\mathcal{F}^{-1}(\varphi(2^{-j}\cdot)\hat{u}(\cdot))={2^{3j}}\int_{\mathbb{R}^3}\check{\varphi}(2^j(x-y))u(y)dy, \ \forall j\in\mathbb{Z},\\
&\dot{S}_ju:=\mathcal{F}^{-1}(\psi(2^{-j}\cdot)\hat{u}(\cdot))={2^{3j}}\int_{\mathbb{R}^3}\check{\psi}(2^j(x-y))u(y)dy, \ \forall j\in\mathbb{Z}.
\end{align*}
It is convenient to introduce the notation $\widetilde{\dot{\Delta }}_ju:=\sum\limits_{|j'-j|\leq2}\dot{\Delta }_{j'}u$, which will use frequently throughout this paper.
\par
In the homogeneous dyadic blocks $\dot{\Delta}_j$, the following Littlewood-Paley decomposition makes sense
\begin{align*}
u=\sum\limits_{j\in\mathbb{Z}}\dot{\Delta}_ju \ \ for\ \ any\ \ u\in \mathcal{S}'_h,
\end{align*}
where $\mathcal{S}'_h$ is given by
\begin{align*}
\mathcal{S}'_h:=\{u\in \mathcal{S}': \lim\limits_{j\rightarrow-\infty}\dot{S}_ju=0\}.
\end{align*}
Moreover, it holds that
\begin{align*}
\dot{S}_ju(x)=\sum\limits_{j'\leq j-1}\dot{\Delta}_{j'}u,\quad j\in \mathbb{Z},\quad  u\in \mathcal{S}'_h.
\end{align*}

We turn to the definition of the Sobolev and Besov space and norms which will come into play in our paper.
\begin{definition}
The homogeneous Sobolev space $\dot{H}^s(\mathbb{R}^3)$ with $s\in \mathbb{R}$ is the space of tempered distributions $u$ over ${\mathbb{R}^3}$, the Fourier transform of which belongs to $L^1_{loc}(\mathbb{R}^3)$ and satisfies
\begin{align*}
\|u\|^2_{H^s}=\int_{\mathbb{R}^d}|\xi|^{2s}|\hat{u}(\xi)|^2d\xi<\infty.
\end{align*}

\end{definition}

\begin{definition}
Let $s\in\mathbb{R}$ and $1\leq p,q\leq\infty$. The homogeneous Besov space $\dot{B}^s_{p,q}$ consists of all tempered distribution $f$ such that
\begin{align*}
\dot{B}^s_{p,q}=\{f\in\mathcal{S}'_h:\|f\|_{\dot{B}^s_{p,q}}<\infty\},
\end{align*}
where
\begin{align*}
		\begin{split}
			\| f\|_{\dot{B}^s_{p,q}}= \left \{
			\begin{array}{ll}
				\left(\sum\limits_{j\in\mathbb{Z}}2^{sjq}\|\dot{\Delta}_jf\|_{L^p}^q\right)^{\frac{1}{q}}                   & if \ 1\leq q <\infty,\\
				\sup\limits_{j\in\mathbb{Z}}2^{sj}\|\dot{\Delta}_jf\|_{L^p}                                & if \ q=\infty.
			\end{array}
			\right.
		\end{split}
	\end{align*}

\end{definition}

As for the case of Lebesgue norms, we have the Bernstein inequalities:
\begin{lemma}\cite{Bahouri} \label{lem}
Let $\mathcal{C}$ be an annulus and $B$ a ball. A constant $C$ exists such that for any nonnegative integer $k$, any couple $(p,q)$ in $[1,\infty]^2$ with $1\leq p\leq q$, and any function $u$ of $L^p$, we have
\begin{align*}
Supp\hat{u}\subset \lambda B&\Rightarrow \|D^ku\|_{L^q}:=\sup\limits_{|\alpha|=k}\|\partial^\alpha u\|_{L^q}\leq C^{k+1}{\lambda^{k+3(\frac{1}{p}-\frac{1}{q})}}\|u\|_{L^p},\\
Supp\hat{u}\subset \lambda \mathcal{C}&\Rightarrow C^{-k-1}\lambda^k\|u\|_{L^p} \leq \|D^ku\|_{L^p}\leq C^{k+1}\lambda^{k}\|u\|_{L^p}.
\end{align*}
\end{lemma}

In the usual manner, we will use $C$ to denote the absolute constants which may be different from line to line unless otherwise stated. Given a Banach space $X$, we denote its norm by $\|\cdot\|_X$. If no otherwise specified, we denote $L^p(\mathbb{R}^3)$ by $L^p$ for convenience. Throughout, we will use the Einstein summation convention and sum over repeated indices.

\section{The proof of Theorem \ref{th1}}
We split $u$ into low and high frequencies: $u=\dot{S}_ku+u^k$. Obviously, there exists a constant $C$ such that
\begin{align*}
\forall p\in[1,\infty], \quad\|\dot{S}_ku\|_{L^p}\leq C\|u\|_{L^p}\quad and\quad \|u^k\|_{L^p}\leq C\|u\|_{L^p}.
\end{align*}
For $u^k$, by the Bernstein inequalities, we have
\begin{align}\label{t1}
\|u^k\|_{L^3}
&\leq\sum\limits_{l\geq k}\|\dot{\Delta}_lu\|_{L^3}
\leq\sum\limits_{l\geq k}2^{\frac{ l}{2}}\|\dot{\Delta}_lu\|_{L^2}\\
&\leq\sum\limits_{l\geq k}2^{-\frac{ l}{2}}\|\nabla \dot{\Delta}_lu\|_{L^2}
\leq 2^{-\frac{ k}{2}}\|\nabla u\|_{L^2}.\notag
\end{align}
This means that the integral $\int_{\mathbb{R}^3}div(\rho u\otimes u)\cdot u^kdx$  is well defined. Indeed, from \eqref{t1}, we obtain
\begin{align}\label{t2}
\int_{\mathbb{R}^3}div(\rho u\otimes u)\cdot u^kdx\leq \|\rho\|_{L^\infty}\|u\|_{L^6}\|\nabla u\|_{L^2}\|u^k\|_{L^3}\leq C(k)\|\rho\|_{L^\infty}\|\nabla u\|_{L^2}^3,
\end{align}
where we used the fact that $div(\rho u)=0$ and $\|u\|_{L^6(\mathbb{R}^3)}\leq C\|\nabla u\|_{L^2(\mathbb{R}^3)}$.

\par
Next, we deduce the $L^{\frac{3}{2}}$ bound for $\nabla P$. Take divergence operator on \eqref{NS1} implies
$$-\Delta P= div(\rho u\cdot\nabla u).$$
It follows that $\nabla P=\nabla(-\Delta)^{-1}div(\rho u\cdot\nabla u)$, which yields
\begin{align}\label{t23}
\|\nabla P\|_{L^{\frac{3}{2}}}\leq C\|\rho u\cdot\nabla u \|_{L^{\frac{3}{2}}}\leq C\|\rho\|_{L^\infty}\|u\|_{L^6}\|\nabla u\|_{L^2}\leq C\|\rho\|_{L^\infty}\|\nabla u\|_{L^2}^2.
\end{align}
Base on \eqref{t1}, \eqref{t2} and \eqref{t23}, taking $L^2$ inner product to \eqref{NS1} with $u^k$, we arrive at
\begin{align}\label{t24}
\int_{\mathbb{R}^3}|\nabla u^k|^2dx=-\int_{\mathbb{R}^3}div(\rho u\otimes u)\cdot u^kdx-\int_{\mathbb{R}^3}\nabla \dot{S}_ku:\nabla u^k dx,
\end{align}
where
\begin{align*}
\int_{\mathbb{R}^3}\nabla P\cdot u^kdx=0.
\end{align*}
Taking $k\rightarrow-\infty$, it is clear that
\begin{align*}
&\lim\limits_{k\rightarrow-\infty}\int_{\mathbb{R}^3}|\nabla u^k|^2dx=\int_{\mathbb{R}^d}|\nabla u|^2dx,\\
&\lim\limits_{k\rightarrow-\infty}|\int_{\mathbb{R}^3} \nabla \dot{S}_ku:\nabla u^k dx|
\leq\|\nabla u\|_{L^2}\lim\limits_{k\rightarrow-\infty}\|\nabla \dot{S}_ku\|_{L^2}=0.
\end{align*}
Therefore, the term \eqref{t24} yields that
\begin{align}\label{t25}
\int_{\mathbb{R}^3}|\nabla u|^2dx&=-\liminf\limits_{k\rightarrow-\infty}\int_{\mathbb{R}^3}div(\rho u\otimes u)\cdot u^kdx=-\liminf\limits_{k\rightarrow-\infty}\int_{\mathbb{R}^3}\rho u\cdot\nabla u\cdot u^kdx.
\end{align}

We will localize the integral in \eqref{t25} into the region near the origin in the frequency space. The first key step is to localize the interactions between $\rho u$ and $u^k$. We split this integral as follows:
\begin{align}\label{t6}
&-\int_{\mathbb{R}^3}\rho u\cdot\nabla u\cdot u^kdx\\
=&-\int_{\mathbb{R}^3}(\rho u_j)\partial_ju_i u^k_idx\notag\\
=&-\int_{\mathbb{R}^3}\dot{S}_k(\rho u_j)\partial_j\dot{S}_ku_i u^k_idx-\int_{\mathbb{R}^3}(\rho u_j)^k\partial_j\dot{S}_ku_i u^k_idx\notag\\
:=&I_1+I_2.\notag
\end{align}
where we used the fact that
\begin{align*}
\int_{\mathbb{R}^3}(\rho u_j)\partial_j u^k_i u_i^k dx=0.
\end{align*}
Next, we will make further efforts to localize $I_1$ and $I_2$  in frequency spaces as follows.
Defining $\widetilde{\dot{\Delta }}_lu:=\sum\limits_{|l'-l|\leq2}\dot{\Delta }_{l'}u$, we rewrite $I_1$ by Bony's decomposition to obtain
\begin{align}\label{t7}
I_1&
=-\int_{\mathbb{R}^3}\dot{S}_k(\rho u_j)\partial_j \dot{S}_k u_i u^k_i dx\\
&=-\int_{\mathbb{R}^3}\partial_j\dot{S}_k u_i\Big(\sum\limits_{l\in\mathbb{Z}}\dot{S}_{l-2}\dot{S}_k (\rho u_j)\dot{\Delta }_l u^k_i     +\sum\limits_{l\in\mathbb{Z}}\dot{\Delta }_l\dot{S}_k (\rho u_j)\dot{S}_{l-2} u^k_i
+\sum\limits_{l\in\mathbb{Z}}\dot{\Delta }_l\dot{S}_k (\rho u_j)\widetilde{\dot{\Delta}}_l u^k_i
\Big)dx\notag\\
&:=I_{11}+I_{12}+I_{13},\notag
\end{align}
For the term $I_{11}$, we note that
\begin{align*}
&\dot{\Delta }_l u^k_i\neq0\Rightarrow l\geq k-1,\\
&supp\widehat{\dot{S}_{l-2}\dot{S}_k (\rho u_j)\dot{\Delta }_l u^k_i}\subset\{\xi:2^{l-2}\leq|\xi|<\frac{9}{8}2^{l+1}\},\\
&supp\widehat{\partial_j\dot{S}_k u_i}\subset\{\xi:|\xi|<2^k\}.
\end{align*}
By virtue of propositions of Littlewood-Paley decomposition, we have
\begin{align*}
I_{11}&=-\int_{\mathbb{R}^3}\partial_j\dot{S}_k u_i\sum\limits_{l=k-1}^{k+1}\dot{S}_{l-2}\dot{S}_k (\rho u_j)\dot{\Delta }_l u^k_idx\\
&=-\int_{\mathbb{R}^3}\sum\limits_{l=k-1}^{k+1}\sum\limits_{l'=l-2}^{k-1}\partial_j\dot{\Delta }_{l'} u_i\dot{S}_{l-2}\dot{S}_k (\rho u_j)\dot{\Delta }_l u^k_idx.
\end{align*}
For the term $I_{12}$, notice that $\dot{\Delta }_l\dot{S}_k (\rho u_j)\neq0\Rightarrow l\leq k$ and $\dot{S}_{l-2} u^k_i\neq0\Rightarrow l\geq k+2$, it means that $\dot{\Delta }_l\dot{S}_k (\rho u_j)\dot{S}_{l-2} u^k_i=0$ for all $l\in\mathbb{Z}$, which yields that
\begin{align*}
I_{12}= -\int_{\mathbb{R}^3}\partial_j\dot{S}_k u_i \sum\limits_{l\in\mathbb{Z}}\dot{\Delta }_l\dot{S}_k (\rho u_j)\dot{S}_{l-2} u^k_idx=0.
\end{align*}
For the term $I_{13}$,  we note that
\begin{align*}
&\dot{\Delta }_l\dot{S}_k (\rho u_j)\neq0 \Rightarrow l\leq k,\\
&\widetilde{\dot{\Delta}}_l u^k_i\neq0 \Rightarrow l\geq k-3,\\
&supp\widehat{\dot{\Delta }_l\dot{S}_k (\rho u_j)\widetilde{\dot{\Delta}}_l u^k_i}\subset\{\xi:|\xi|<5\times2^{l+1}\}.
\end{align*}
It follows that
\begin{align*}
I_{13}=-\int_{\mathbb{R}^d}\partial_j\dot{S}_k u_i \sum\limits_{l=k-3}^k\dot{\Delta }_l\dot{S}_k (\rho u_j)\widetilde{\dot{\Delta}}_l u^k_i dx.
\end{align*}
Base on above estimates, \eqref{t7} yields
\begin{align}\label{t8}
&I_1\\=&-\int_{\mathbb{R}^3}\sum\limits_{l=k-1}^{k+1}\sum\limits_{l'=l-2}^{k-1}\partial_j\dot{\Delta }_{l'} u_i\dot{S}_{l-2}\dot{S}_k (\rho u_j)\dot{\Delta }_l u^k_idx-\int_{\mathbb{R}^d}\partial_j\dot{S}_k u_i \sum\limits_{l=k-3}^k\dot{\Delta }_l\dot{S}_k (\rho u_j)\widetilde{\dot{\Delta}}_l u^k_i dx.\notag
\end{align}
We rewrite $I_2$ by Bony's decomposition as
\begin{align}\label{t9}
I_2&
=-\int_{\mathbb{R}^3}\partial_j\dot{S}_k u_i(\rho u_j)^ku^k_idx\\
&=-\int_{\mathbb{R}^3}\partial_j\dot{S}_k u_i(\sum\limits_{l\in\mathbb{Z}}\dot{S}_{l-2}(\rho u_j)^k\dot{\Delta }_lu^k_i
+\sum\limits_{l\in\mathbb{Z}}\dot{S}_{l-2}u^k_i\dot{\Delta }_l(\rho u_j)^k
+\sum\limits_{l\in\mathbb{Z}}\dot{\Delta }_l(\rho u_j)^k\widetilde{\dot{\Delta}}_lu^k_i)dx\notag\\
&:=I_{21}+I_{22}+I_{23}.\notag
\end{align}
We estimate $I_{21}$, $I_{22}$ and $I_{23}$ as follows.
Notice that
\begin{align*}
&\dot{S}_{l-2}(\rho u_j)^k\neq0\Rightarrow l\geq k+2,\\
&\dot{\Delta }_lu^k_i\neq0\Rightarrow l\geq k-1,\\
&supp\widehat{\dot{S}_{l-2}u^k_j\dot{\Delta }_lu^k_i}\subset\{\xi:2^{l-2}\leq|\xi|<\frac{9}{8}2^{l+1}\}\quad for\quad l\geq k+2,\\
&supp\widehat{\partial_j\dot{S}_k u_i}\subset\{\xi:|\xi|<2^k\},
\end{align*} it follows that
\begin{align*}
supp\widehat{\dot{S}_{l-2}u^k_j\dot{\Delta }_lu^k_i}\cap supp\widehat{\partial_j\dot{S}_k u_i}=\emptyset\quad for \quad l\geq k+2.
\end{align*}
Therefore we have
\begin{align*}
I_{21}=-\int_{\mathbb{R}^3}\partial_j\dot{S}_k u_i\sum\limits_{l\geq k+2}\dot{S}_{l-2}(\rho u_j)^k\dot{\Delta }_lu^k_idx=0.
\end{align*}
A similar arguments also hold for $I_{22}$, which shows that
\begin{align*}
I_{22}=0.
\end{align*}
\par
For the term $I_{23}$, since
\begin{align*}
&\dot{\Delta }_l (\rho u_j)^k\neq0 \Rightarrow l\geq k-1,\\
&\widetilde{\dot{\Delta}}_l u_i^k\neq0 \Rightarrow l\geq k-3,\\
&supp\widehat{\dot{\Delta }_l (\rho u_j)^k\widetilde{\dot{\Delta}}_l u^k_i}\subset\{\xi:|\xi|<5\times2^{l+1}\},
\end{align*}
then we have
\begin{align*}
I_{23}=-\int_{\mathbb{R}^3}\partial_j\dot{S}_k u_i\sum\limits_{l\geq k-1}\dot{\Delta }_l(\rho u_j)^k\widetilde{\dot{\Delta}}_lu^k_idx.
\end{align*} It follows that
\begin{align}\label{+1}
I_2=-\int_{\mathbb{R}^3}\partial_j\dot{S}_k u_i\sum\limits_{l\geq k-1}\dot{\Delta }_l(\rho u_j)^k\widetilde{\dot{\Delta}}_lu^k_idx.
\end{align}
Combining \eqref{t8}, \eqref{+1} and \eqref{t25} yields
\begin{align}\label{t10}
&\int_{\mathbb{R}^3}|\nabla u|^2dx\\
=&{\color{red}-}\liminf\limits_{k\rightarrow-\infty}\Big\{ \int_{\mathbb{R}^3}\sum\limits_{l=k-1}^{k+1}\sum\limits_{l'=l-2}^{k-1}\partial_j\dot{\Delta }_{l'} u_i\dot{S}_{l-2}\dot{S}_k (\rho u_j)\dot{\Delta }_l u^k_idx
+\int_{\mathbb{R}^3}\partial_j\dot{S}_k u_i \sum\limits_{l=k-3}^k\dot{\Delta }_l\dot{S}_k (\rho u_j)\widetilde{\dot{\Delta}}_l u^k_i dx\notag\\
&~~~~~~~~~~~+\int_{\mathbb{R}^3}\partial_j\dot{S}_k u_i\sum\limits_{l\geq k-1}\dot{\Delta }_l(\rho u_j)^k\widetilde{\dot{\Delta}}_lu^k_idx  \Big\}.\notag
\end{align}
In view of the formula \eqref{t10}, the second key step is to localize the interactions between the density $\rho$ and the velocity $u$ and give suitable estimates on different frequency parts.
We achieve our goals through the following three propositions.
\begin{proposition}\label{pro1}
Assume $(\rho,u)$ be a smooth solution of \eqref{NS1} satisfying the integrability conditions in Theorem \ref{th1}, we have
\begin{align*}
\liminf\limits_{k\rightarrow-\infty}\Big\{ \int_{\mathbb{R}^d}\sum\limits_{l=k-1}^{k+1}\sum\limits_{l'=l-2}^{k-1}\partial_j\dot{\Delta }_{l'} u_i\dot{S}_{l-2}\dot{S}_k (\rho u_j)\dot{\Delta }_l u^k_idx\Big\}=0.
\end{align*}
\end{proposition}
\emph{Proof.} By Bony's decomposition, we get
\begin{align*}
&\int_{\mathbb{R}^3}\sum\limits_{l=k-1}^{k+1}\sum\limits_{l'=l-2}^{k-1}\partial_j\dot{\Delta }_{l'} u_i\dot{S}_{l-2}\dot{S}_k (\rho u_j)\dot{\Delta }_l u^k_idx\\
=&\int_{\mathbb{R}^d}\sum\limits_{l=k-1}^{k+1}\sum\limits_{l'=l-2}^{k-1}\partial_j\dot{\Delta }_{l'} u_i\dot{S}_{l-2}\dot{S}_k \Big(
\sum\limits_{l''\in\mathbb{Z}}\dot{S}_{l''-2}\rho\dot{\Delta }_{l''} u_j+\sum\limits_{l''\in\mathbb{Z}}\dot{\Delta }_{l''}\rho\dot{S}_{l''-2} u_j
+\sum\limits_{l''\in\mathbb{Z}}\widetilde{\dot{\Delta }}_{l''}\rho\dot{\Delta }_{l''} u_j\Big)\dot{\Delta }_l u^k_idx\notag\\
:=&J_1+J_2+J_3.\notag
\end{align*}
For the term $J_1$, we observe that
\begin{align*}
&supp\widehat{\dot{S}_{l''-2}\rho\dot{\Delta }_{l''} u_j}\subset\{\xi:2^{l''-2}\leq|\xi|<\frac{9}{8}2^{l''+1}\},\\
&\dot{S}_{l-2}\dot{S}_k (\dot{S}_{l''-2}\rho\dot{\Delta }_{l''} u_j)\neq0\Rightarrow l''\leq k+1.
\end{align*}
Hence it follows from the above facts that
\begin{align*}
J_1=&\int_{\mathbb{R}^3}\sum\limits_{l=k-1}^{k+1}\sum\limits_{l'=l-2}^{k-1}\sum\limits_{l''\leq k+1}\partial_j\dot{\Delta }_{l'} u_i\dot{S}_{l-2}\dot{S}_k (\dot{S}_{l''-2}\rho\dot{\Delta }_{l''} u_j)\dot{\Delta }_l u^k_idx\\
\leq&\sum\limits_{l=k-1}^{k+1}\sum\limits_{l'=l-2}^{k-1}\sum\limits_{l''\leq k+1}\|\dot{S}_{l''-2}\rho\|_{L^\infty}\|\dot{\Delta }_{l''} u\|_{L^\infty}\|\nabla\dot{\Delta }_{l'} u\|_{L^2}\|\dot{\Delta }_l u\|_{L^2}\\
\leq&C\|\rho\|_{L^\infty}\sum\limits_{l''\leq k+1}2^{l''-k}(2^{-l''}\|\dot{\Delta }_{l''} u\|_{L^\infty})\sum\limits_{l=k-1}^{k+1}\sum\limits_{l'=l-2}^{k-1}\|\nabla\dot{\Delta }_{l'} u\|_{L^2}2^{k-l}2^l\|\dot{\Delta }_l u\|_{L^2}\\
\leq&C\|\rho\|_{L^\infty}\|\dot{S}_{k} u\|_{\dot{B}^{-1}_{\infty,\infty}}\|\nabla\dot{S}_{k+2}u\|_{L^2}^2,
\end{align*}where we have applied Young's inequality and H\"{o}lder's inequality.
Applying the similar arguments, we estimate $J_2$ as follows
\begin{align*}
J_{2}=&\int_{\mathbb{R}^3}\sum\limits_{l=k-1}^{k+1}\sum\limits_{l'=l-2}^{k-1}\sum\limits_{l''\leq k+1 }\partial_j\dot{\Delta }_{l'} u_i\dot{S}_{l-2}\dot{S}_k (\dot{\Delta }_{l''}\rho\dot{S}_{l''-2} u_j)\dot{\Delta }_l u^k_idx\\
\leq& \sum\limits_{l=k-1}^{k+1}\sum\limits_{l'=l-2}^{k-1}\sum\limits_{ l''\leq k+1 }\|\dot{\Delta }_{l''}\rho\|_{L^\infty}   \|\dot{S}_{l''-2} u\|_{L^\infty} \|\nabla \dot{\Delta }_{l'}u\|_{L^2}\|\dot{\Delta }_l u\|_{L^2}\\
\leq&C\|\rho\|_{L^\infty} \sum\limits_{l''\leq k+1 }2^{l''-k}\big(2^{-(l''-2)} \|\dot{S}_{l''-2} u\|_{L^\infty}\big)\sum\limits_{l=k-1}^{k+1}\sum\limits_{l'=l-2}^{k-1}\|\nabla \dot{\Delta }_{l'}u\|_{L^2}2^{k-l}2^l\|\dot{\Delta }_l u\|_{L^2}\\
 \leq&C\|\rho\|_{L^\infty}
\|\dot{S}_{k} u\|_{\dot{B}^{-1}_{\infty,\infty}}\|\nabla\dot{S}_{k+2}u\|_{L^2}^2,
\end{align*} where we have used the fact that
\begin{align}\label{t11}
2^{-k}\|\dot{S}_ku\|_{L^\infty}&\leq 2^{-k}\sum\limits_{l\leq k-1}\|\dot{\Delta}_lu\|_{L^\infty}\\
&\leq\sum\limits_{l\leq k-1}2^{l-k}2^{-l}\|\dot{\Delta}_lu\|_{L^\infty}\notag\\
&\leq C\|\dot{S}_ku\|_{\dot{B}_{\infty,\infty}^{-1}}.\notag
\end{align}
For the term $J_3$, we observe that
\begin{align*}
&supp\widehat{\widetilde{\dot{\Delta }}_{l''}\rho\dot{\Delta }_{l''} u_j}\subset\{\xi:|\xi|<5\times2^{l''+1}\},\\
&\dot{S}_{l-2}\dot{S}_k (\widetilde{\dot{\Delta }}_{l''}\rho\dot{\Delta }_{l''} u_j)\neq 0 \quad for~ all~l''\in\mathbb{Z}.
\end{align*}
Thus, we decompose $J_3$ as follows.
\begin{align*}
J_3=&\int_{\mathbb{R}^3}\sum\limits_{l=k-1}^{k+1}\sum\limits_{l'=l-2}^{k-1}\sum\limits_{l''\in\mathbb{Z}}\partial_j\dot{\Delta }_{l'} u_i\dot{S}_{l-2}\dot{S}_k (\widetilde{\dot{\Delta }}_{l''}\rho\dot{\Delta }_{l''} u_j)\dot{\Delta }_l u^k_idx\\
=& \int_{\mathbb{R}^3}\sum\limits_{l=k-1}^{k+1}\sum\limits_{l'=l-2}^{k-1}\sum\limits_{l''\leq k}\partial_j\dot{\Delta }_{l'} u_i\dot{S}_{l-2}\dot{S}_k (\widetilde{\dot{\Delta }}_{l''}\rho\dot{\Delta }_{l''} u_j)\dot{\Delta }_l u^k_idx\\
&+ \int_{\mathbb{R}^3}\sum\limits_{l=k-1}^{k+1}\sum\limits_{l'=l-2}^{k-1}\sum\limits_{l''\geq k+1}\partial_j\dot{\Delta }_{l'} u_i\dot{S}_{l-2}\dot{S}_k (\widetilde{\dot{\Delta }}_{l''}\rho\dot{\Delta }_{l''} u_j)\dot{\Delta }_l u^k_idx\\
:=&J_{31}+J_{32}.
\end{align*}
We estimate $J_{31}$ and $J_{32}$ as follows.
\begin{align*}
J_{31}
&\leq\sum\limits_{l=k-1}^{k+1}\sum\limits_{l'=l-2}^{k-1}\sum\limits_{l''\leq k}\|\widetilde{\dot{\Delta }}_{l''}\rho\|_{L^\infty}\|\dot{\Delta }_{l''} u\|_{L^\infty}\|\nabla\dot{\Delta }_{l'} u\|_{L^2}  \|\dot{\Delta }_l u\|_{L^2}\\
&\leq C\|\rho\|_{L^\infty}\sum\limits_{l''\leq k}2^{l''-k}2^{-l''}\|\dot{\Delta }_{l''} u\|_{L^\infty}\sum\limits_{l=k-1}^{k+1}\sum\limits_{l'=l-2}^{k-1}\|\nabla\dot{\Delta }_{l'} u\|_{L^2}2^{k-l} 2^l \|\dot{\Delta }_l u\|_{L^2}\\
&\leq C\|\rho\|_{L^\infty}
\|\dot{S}_{k+1} u\|_{\dot{B}_{\infty,\infty}^{-1}}\|\nabla\dot{S}_{k+2}u\|_{L^2}^2.
\end{align*}
\begin{align*}
J_{32}\leq&\sum\limits_{l=k-1}^{k+1}\sum\limits_{l'=l-2}^{k-1}\sum\limits_{l''\geq k+1}\|\widetilde{\dot{\Delta }}_{l''}\rho\|_{L^\infty}\|\dot{\Delta }_{l''} u\|_{L^2}\|\nabla\dot{\Delta }_{l'} u\|_{L^2}   \|\dot{\Delta }_l u\|_{L^\infty}\\
\leq&C\|\rho\|_{L^\infty}\sum\limits_{l''\geq k+1}2^{k-l''}2^{l''}\|\dot{\Delta }_{l''} u\|_{L^2}\sum\limits_{l=k-1}^{k+1}\sum\limits_{l'=l-2}^{k-1}\|\nabla\dot{\Delta }_{l'} u\|_{L^2}  2^{l-k}\big(2^{-l}\|\dot{\Delta }_l u\|_{L^\infty}\big)\\
\leq&C\|\rho\|_{L^\infty}\|\nabla u\|_{L^2}\|\dot{S}_{k+2} u\|_{\dot{B}_{\infty,\infty}^{-1}}\sum\limits_{l=k-1}^{k+1}\sum\limits_{l'=l-2}^{k-1}2^{l-k}\|\nabla\dot{\Delta }_{l'} u\|_{L^2}\\
\leq&C\|\rho\|_{L^\infty}\|\nabla u\|_{L^2}\|\dot{S}_{k+2} u\|_{\dot{B}_{\infty,\infty}^{-1}}\|\nabla\dot{S}_{k}u\|_{L^2}.
\end{align*}
Combining the condition \eqref{con1}, $\|\nabla u\|_{L^2}<\infty$ and $\lim\limits_{k\rightarrow-\infty}\|\nabla\dot{S}_{k}u\|_{L^2}=0$,
we conclude that
$$\liminf\limits_{k\rightarrow-\infty} \{J_1+J_2+J_{3}\}=0.$$
This completes the proof of Proposition \ref{pro1}.\qquad $\hfill\Box$\\

\begin{proposition}\label{pro2}
Assume $(\rho,u)$ be a smooth solution of \eqref{NS1} satisfying the integrability conditions in Theorem \ref{th1}, we have
\begin{align*}
\liminf\limits_{k\rightarrow-\infty}\Big\{ \int_{\mathbb{R}^3}\partial_j\dot{S}_k u_i \sum\limits_{l=k-3}^k\dot{\Delta }_l\dot{S}_k (\rho u_j)\widetilde{\dot{\Delta}}_l u^k_i dx\Big\}=0.
\end{align*}
\end{proposition}
\emph{Proof.} By Bony's decomposition, we have
\begin{align*}
&\int_{\mathbb{R}^3}\partial_j\dot{S}_k u_i \sum\limits_{l=k-3}^k\dot{\Delta }_l\dot{S}_k (\rho u_j)\widetilde{\dot{\Delta}}_l u^k_i dx\\
=&\int_{\mathbb{R}^3}\partial_j\dot{S}_k u_i \sum\limits_{l=k-3}^k\dot{\Delta }_l\dot{S}_k \Big(
\sum\limits_{l'\in\mathbb{Z}}\dot{S}_{l'-2}\rho\dot{\Delta }_{l'} u_j+\sum\limits_{l'\in\mathbb{Z}}\dot{\Delta }_{l'}\rho\dot{S}_{l'-2} u_j
+\sum\limits_{l'\in\mathbb{Z}}\widetilde{\dot{\Delta }}_{l'}\rho\dot{\Delta }_{l'} u_j\Big)\widetilde{\dot{\Delta}}_l u^k_i dx\\
:=&K_1+K_2+K_3.
\end{align*}
For the term $K_1$, we observe that
\begin{align*}
&supp\widehat{\dot{S}_{l'-2}\rho\dot{\Delta }_{l'} u_j}\subset\{\xi:2^{l'-2}\leq|\xi|<\frac{9}{8}2^{l'+1}\},\\
&\dot{\Delta }_l\dot{S}_k (\dot{S}_{l'-2}\rho\dot{\Delta }_{l'} u_j)\neq0\Rightarrow l'\leq k+1~and~|l'-l|\leq4.
\end{align*}
Therefore, from H\"{o}lder's inequality, Young's inequality and the fact \eqref{t11}, it follows that
\begin{align*}
K_1=&\int_{\mathbb{R}^3}\sum\limits_{l=k-3}^k\sum\limits_{l'=k-7}^{k+1}\partial_j\dot{S}_k u_i\dot{\Delta }_l\dot{S}_k(\dot{S}_{l'-2}\rho\dot{\Delta }_{l'} u_j)\widetilde{\dot{\Delta}}_l u^k_i dx\\
\leq&\sum\limits_{l=k-3}^k\sum\limits_{l'=k-7}^{k+1}\|\dot{S}_{l'-2}\rho\|_{L^\infty}\|\dot{\Delta }_{l'} u\|_{L^2}\|\nabla\dot{S}_k u\|_{L^\infty}\|\widetilde{\dot{\Delta}}_l u\|_{L^2}\\
\leq& C \|\rho\|_{L^\infty} 2^{-k}\|\dot{S}_k u\|_{L^\infty}\sum\limits_{l=k-3}^k\sum\limits_{l'=k-7}^{k+1}2^{k-l'}2^{l'}\|\dot{\Delta }_{l'} u\|_{L^2}2^{k-l}2^l\|\widetilde{\dot{\Delta}}_l u\|_{L^2}\\
\leq& C \|\rho\|_{L^\infty} 2^{-k}\|\dot{S}_k u\|_{L^\infty}\|\nabla\dot{S}_{k+2}u\|_{L^2}^2\\
\leq& C \|\rho\|_{L^\infty} \|\dot{S}_ku\|_{\dot{B}_{\infty,\infty}^{-1}}\|\nabla\dot{S}_{k+2}u\|_{L^2}^2,
\end{align*}

Similarly, we estimate $K_2$ as follows.
\begin{align*}
K_2=&\int_{\mathbb{R}^3}\sum\limits_{l=k-3}^k\sum\limits_{l'=k-7}^{k+1}\partial_j\dot{S}_k u_i\dot{\Delta }_l\dot{S}_k(\dot{\Delta }_{l'}\rho\dot{S}_{l'-2} u_j)\widetilde{\dot{\Delta}}_l u^k_i dx\\
\leq&C\sum\limits_{l=k-3}^k\sum\limits_{l'=k-7}^{k+1}\|\nabla\dot{S}_k u\|_{L^2}\|\dot{\Delta }_{l'}\rho\|_{L^\infty}\|\dot{S}_{l'-2} u\|_{L^\infty}\|\widetilde{\dot{\Delta}}_l u\|_{L^2}\\
\leq& C \|\rho\|_{L^\infty}\|\nabla\dot{S}_k u\|_{L^2} \sum\limits_{l'=k-7}^{k+1}2^{l'-k}2^{-(l'-2)}\|\dot{S}_{l'-2} u\|_{L^\infty}\sum\limits_{l=k-3}^k2^{k-l}2^l\|\widetilde{\dot{\Delta}}_l u\|_{L^2}\\
\leq& C \|\rho\|_{L^\infty} \|\dot{S}_{k} u\|_{\dot{B}_{\infty,\infty}^{-1}}\|\nabla\dot{S}_{k+1} u\|_{L^2}^2.
\end{align*}

For the term $K_3$, we observe that
\begin{align*}
&supp\widehat{\widetilde{\dot{\Delta }}_{l'}\rho\dot{\Delta }_{l'} u_j}\subset\{\xi:|\xi|<5\times2^{l'+1}\},\\
&\dot{\Delta }_l\dot{S}_k (\widetilde{\dot{\Delta }}_{l'}\rho\dot{\Delta }_{l'} u_j)\neq0\Rightarrow l'\geq l-3.
\end{align*}
From direct calculations, we show that
\begin{align*}
K_3=&\int_{\mathbb{R}^3}\sum\limits_{l=k-3}^k\sum\limits_{l'\geq l-3}\partial_j\dot{S}_k u_i\dot{\Delta }_l\dot{S}_k(\widetilde{\dot{\Delta }}_{l'}\rho\dot{\Delta }_{l'} u_j)\widetilde{\dot{\Delta}}_l u^k_i dx\\
\leq&C2^k\|\dot{S}_ku\|_{L^\infty} \sum\limits_{l=k-3}^k\sum\limits_{l'\geq l-3}\|\dot{\Delta }_{l'}\rho\|_{L^\infty}\|\dot{\Delta }_{l'} u\|_{L^2}\|\widetilde{\dot{\Delta}}_l u\|_{L^2}\\
\leq&C2^{-k}\|\dot{S}_ku\|_{L^\infty} \|\rho\|_{L^\infty}\sum\limits_{l=k-3}^k\sum\limits_{l'\geq l-3}2^{l-l'}2^{l'}\|\dot{\Delta }_{l'} u\|_{L^2}2^{2(k-l)}2^{l}\|\widetilde{\dot{\Delta}}_l u\|_{L^2}\\
\leq&C2^{-k}\|\dot{S}_ku\|_{L^\infty} \|\rho\|_{L^\infty}\|\nabla u\|_{L^2}\sum\limits_{l=k-3}^k2^{2(k-l)}2^{l}\|\widetilde{\dot{\Delta}}_l u\|_{L^2}\\
\leq&C \|\rho\|_{L^\infty}\|\dot{S}_{k} u\|_{\dot{B}_{\infty,\infty}^{-1}}\|\nabla u\|_{L^2}\|\nabla\dot{S}_{k+1} u\|_{L^2}.
\end{align*}
Combining the conditions \eqref{con1}, $\|\nabla u\|_{L^2}<\infty$ and $\lim\limits_{k\rightarrow-\infty}\|\nabla\dot{S}_{k}u\|_{L^2}=0$,
we conclude that
$$\liminf\limits_{k\rightarrow-\infty}\{ K_1+K_2+K_3\}=0.$$
This completes the proof of Proposition \ref{pro2}.\qquad $\hfill\Box$\\

\begin{proposition}\label{pro3}
Assume $(\rho,u)$ be a smooth solution of \eqref{NS1} satisfying the integrability conditions in Theorem \ref{th1}, we have
\begin{align*}
\liminf\limits_{k\rightarrow-\infty}\Big\{ \int_{\mathbb{R}^3}\partial_j\dot{S}_k u_i\sum\limits_{l\geq k-1}\dot{\Delta }_l(\rho u_j)^k\widetilde{\dot{\Delta}}_lu^k_idx\Big\}=0.
\end{align*}
\end{proposition}
\emph{Proof.} By Bony's decomposition, we have
\begin{align*}
&\int_{\mathbb{R}^3}\partial_j\dot{S}_k u_i\sum\limits_{l\geq k-1}\dot{\Delta }_l(\rho u_j)^k\widetilde{\dot{\Delta}}_lu^k_idx\\
=&\int_{\mathbb{R}^3}\partial_j\dot{S}_k u_i\sum\limits_{l\geq k-1}\dot{\Delta }_l\Big(
\sum\limits_{l'\in\mathbb{Z}}\dot{S}_{l'-2}\rho\dot{\Delta }_{l'} u_j+\sum\limits_{l'\in\mathbb{Z}}\dot{\Delta }_{l'}\rho\dot{S}_{l'-2} u_j
+\sum\limits_{l'\in\mathbb{Z}}\widetilde{\dot{\Delta }}_{l'}\rho\dot{\Delta }_{l'} u_j\Big)^k\widetilde{\dot{\Delta}}_lu^k_idx\\
:=&L_1+L_2+L_3.
\end{align*}

For the term $L_1$, we observe that
\begin{align*}
&supp\widehat{\dot{S}_{l'-2}\rho\dot{\Delta }_{l'} u_j}\subset\{\xi:2^{l'-2}\leq|\xi|<\frac{9}{8}2^{l'+1}\},\\
&\dot{\Delta }_l(\dot{S}_{l'-2}\rho\dot{\Delta }_{l'} u_j)^k\neq0\Rightarrow l'\geq k-1~and~|l'-l|\leq4.
\end{align*}
Let $\theta\in(0,1)$ whose value will be determined later. We decompose $L_1$ as follows.
\begin{align*}
L_1
=&\int_{\mathbb{R}^3}\sum\limits_{l\geq k-1}\sum\limits_{|l'-l|\leq4}\partial_j\dot{S}_k u_i\dot{\Delta }_l\Big(
\dot{S}_{l'-2}\rho\dot{\Delta }_{l'} u_j\Big)^k\widetilde{\dot{\Delta}}_lu^k_idx\\
=&\int_{\mathbb{R}^3}\sum\limits_{l\geq k-1}\sum\limits_{k-5\leq l'\leq[\theta k]-1}\partial_j\dot{S}_k u_i\dot{\Delta }_l\Big(
\dot{S}_{l'-2}\rho\dot{\Delta }_{l'} u_j\Big)^k\widetilde{\dot{\Delta}}_lu^k_idx\\
&+\int_{\mathbb{R}^3}\sum\limits_{l\geq k-1}\sum\limits_{ l'\geq[\theta k]}\partial_j\dot{S}_k u_i\dot{\Delta }_l\Big(
\dot{S}_{l'-2}\rho\dot{\Delta }_{l'} u_j\Big)^k\widetilde{\dot{\Delta}}_lu^k_idx\\
=&L_{11}+L_{12}.
\end{align*}
For the term $L_{12}$, by using H\"{o}lder's inequality, Young's inequality and Lemma \ref{lem} we have
\begin{align*}
L_{12}\leq&\|\nabla\dot{S}_k u\|_{L^\infty}\sum\limits_{l\geq k-1}\sum\limits_{ l'\geq[\theta k]}\|\dot{S}_{l'-2}\rho\|_{L^\infty}\|\dot{\Delta }_{l'} u\|_{L^2} \|\widetilde{\dot{\Delta}}_lu\|_{L^2}\\
\leq&C2^{\frac{3}{2}k}\|\dot{S}_k u\|_{L^6}\|\rho\|_{L^\infty}\sum\limits_{l\geq k-1}\sum\limits_{ l'\geq[\theta k]}2^{l'}\|\dot{\Delta }_{l'} u\|_{L^2}2^l \|\widetilde{\dot{\Delta}}_lu\|_{L^2}2^{-l-l'}\\
\leq&C2^{(\frac{1}{2}-\theta)k}\|\dot{S}_k u\|_{L^6}\|\rho\|_{L^\infty}\sum\limits_{l\geq k-1}\sum\limits_{ l'\geq[\theta k]}2^{(\theta k-l')}2^{l'}\|\dot{\Delta }_{l'} u\|_{L^2}2^{(k-l)}2^l \|\widetilde{\dot{\Delta}}_lu\|_{L^2}\\
\leq&C2^{(\frac{1}{2}-\theta)k}\|\nabla\dot{S}_k u\|_{L^2}\|\rho\|_{L^\infty}\|\nabla u\|_{L^2}^2,
\end{align*}
where we have used the fact $H^1(\mathbb{R}^3)\hookrightarrow L^6(\mathbb{R}^3)$. Choosing $\theta=\frac{1}{2}$, since
$$\lim\limits_{k\rightarrow-\infty}\|\dot{S}_k \nabla u\|_{L^2}=0,$$
then we deduce
$$\liminf\limits_{k\rightarrow-\infty} L_{12}=0.$$
Then, by similar calculations, the quantity $L_{11}$ can be bounded as follows.
\begin{align*}
L_{11}=&\int_{\mathbb{R}^3}\sum\limits_{l\geq k-1}\sum\limits_{k-5\leq l'\leq[\frac{ k}{2}]-1}\partial_j\dot{S}_k u_i\dot{\Delta }_l\Big(
\dot{S}_{l'-2}\rho\dot{\Delta }_{l'} u_j\Big)^k\widetilde{\dot{\Delta}}_lu^k_idx\\
\leq&C2^{-k}\|\dot{S}_k  u\|_{L^\infty} \|\rho\|_{L^\infty}\sum\limits_{l\geq k-1}2^{k-l}2^l \|\widetilde{\dot{\Delta}}_lu\|_{L^2}\sum\limits_{k-5\leq l'\leq[\frac{ k}{2}]-1}2^{k-l'}2^{l'}\|\dot{\Delta }_{l'} u\|_{L^2}\\
\leq&C2^{-k}\|\dot{S}_k  u\|_{L^\infty} \|\rho\|_{L^\infty} \|\nabla u\|_{L^2}\|\nabla\dot{S}_{[\frac{ k}{2}]} u\|_{L^2}\\
\leq&C \|\rho\|_{L^\infty}\|\dot{S}_{k} u\|_{\dot{B}_{\infty,\infty}^{-1}} \|\nabla u\|_{L^2}\|\nabla\dot{S}_{[\frac{ k}{2}]} u\|_{L^2},
\end{align*}where we also used the fact \eqref{t11}.

We now consider the term $L_2$. For localizing the nonlinear interactions between $\rho$ and $u$ in $L_2$, we decompose it into three parts:
\begin{align*}
L_2=&\int_{\mathbb{R}^3}\sum\limits_{l\geq k-1}\sum\limits_{|l'-l|\leq 4}\partial_j\dot{S}_k u_i\dot{\Delta }_l\Big(
\dot{\Delta }_{l'}\rho\dot{S}_{l'-2} u_j\Big)^k\widetilde{\dot{\Delta}}_lu^k_idx\\
=&\int_{\mathbb{R}^3}\sum\limits_{l\geq k-1}\sum\limits_{|l'-l|\leq 4}\sum\limits_{l''\leq l'-3}\partial_j\dot{S}_k u_i\dot{\Delta }_l\Big(
\dot{\Delta }_{l'}\rho\dot{\Delta}_{l''} u_j\Big)^k\widetilde{\dot{\Delta}}_lu^k_idx\\
\leq&\int_{\mathbb{R}^3}\sum\limits_{l\geq k-1}\sum\limits_{|l'-l|\leq 4}\sum\limits_{l''\leq k-1}\partial_j\dot{S}_k u_i\dot{\Delta }_l\Big(
\dot{\Delta }_{l'}\rho\dot{\Delta}_{l''} u_j\Big)^k\widetilde{\dot{\Delta}}_lu^k_idx\\
&+\int_{\mathbb{R}^3}\sum\limits_{l\geq k-1}\sum\limits_{|l'-l|\leq 4}\sum\limits_{l''\geq [\frac{k}{2}]}\partial_j\dot{S}_k u_i\dot{\Delta }_l\Big(
\dot{\Delta }_{l'}\rho\dot{\Delta}_{l''} u_j\Big)^k\widetilde{\dot{\Delta}}_lu^k_idx\\
&+\int_{\mathbb{R}^3}\sum\limits_{l\geq k-1}\sum\limits_{|l'-l|\leq 4}\sum\limits_{k\leq l''\leq [\frac{k}{2}]}\partial_j\dot{S}_k u_i\dot{\Delta }_l\Big(
\dot{\Delta }_{l'}\rho\dot{\Delta}_{l''} u_j\Big)^k\widetilde{\dot{\Delta}}_lu^k_idx\\
:=&L_{21}+L_{22}+L_{23}.
\end{align*}

Apply Bernstein inequalities, Young inequality and the fact $H^1(\mathbb{R}^3)\hookrightarrow L^6(\mathbb{R}^3)$, we have
\begin{align*}
L_{21}\leq&C\sum\limits_{l\geq k-1}\sum\limits_{|l'-l|\leq 4}\sum\limits_{l''\leq k-1}
\|\nabla\dot{S}_k u\|_{L^2}
\|\dot{\Delta }_{l'}\rho\|_{L^\infty}\|\dot{\Delta}_{l''} u\|_{L^\infty}\|\widetilde{\dot{\Delta}}_lu\|_{L^2}\\
\leq&C\|\rho\|_{L^\infty}\|\nabla\dot{S}_k u\|_{L^2}\sum\limits_{l\geq k-1}\sum\limits_{|l'-l|\leq 4}2^{k-l}2^{l}\|\widetilde{\dot{\Delta}}_lu\|_{L^2}\sum\limits_{l''\leq k-1}
2^{l''-k}2^{-l''}\|\dot{\Delta}_{l''} u\|_{L^\infty}\\
\leq&C\|\rho\|_{L^\infty}\|\nabla\dot{S}_k u\|_{L^2}\|\nabla u\|_{L^2}\|\dot{S}_{k} u\|_{\dot{B}_{\infty,\infty}^{-1}}.
\end{align*}

\begin{align*}
L_{22}\leq&C\|\rho\|_{L^\infty}\|\nabla\dot{S}_k u\|_{L^\infty}
\sum\limits_{l\geq k-1}\sum\limits_{|l'-l|\leq 4}\sum\limits_{l''\geq [\frac{k}{2}]}
\|\dot{\Delta}_{l''} u\|_{L^2}\|\widetilde{\dot{\Delta}}_lu\|_{L^2}\\
\leq&C\|\rho\|_{L^\infty}2^{\frac{3}{2}k}\|\dot{S}_k u\|_{L^6}
\sum\limits_{l\geq k-1}\sum\limits_{|l'-l|\leq 4}\sum\limits_{l''\geq [\frac{k}{2}]}
2^{l''}\|\dot{\Delta}_{l''} u\|_{L^2}2^l\|\widetilde{\dot{\Delta}}_lu\|_{L^2}2^{-l''-l}\\
\leq&C\|\rho\|_{L^\infty}\|\dot{S}_k u\|_{L^6}
\sum\limits_{l\geq k-1}\sum\limits_{|l'-l|\leq 4}\sum\limits_{l''\geq [\frac{k}{2}]}
2^{(\frac{k}{2}-l'')}2^{l''}\|\dot{\Delta}_{l''} u\|_{L^2}2^{(k-l)}2^l\|\widetilde{\dot{\Delta}}_lu\|_{L^2}\\
\leq&C\|\rho\|_{L^\infty}\|\nabla \dot{S}_k u\|_{L^2}\|\nabla u\|_{L^2}^2.
\end{align*}

\begin{align*}
L_{23}\leq&C\|\rho\|_{L^\infty} \|\nabla\dot{S}_k u\|_{L^\infty}  \sum\limits_{l\geq k-1}\sum\limits_{|l'-l|\leq 4}\sum\limits_{k\leq l''\leq [\frac{k}{2}]-1}
\|\dot{\Delta}_{l''} u\|_{L^2}\|\widetilde{\dot{\Delta}}_lu\|_{L^2}\\
\leq&C\|\rho\|_{L^\infty} 2^{-k}\|\dot{S}_k u\|_{L^\infty}  \sum\limits_{l\geq k-1}\sum\limits_{|l'-l|\leq 4}\sum\limits_{k\leq l''\leq [\frac{k}{2}]-1}
2^{k-l''}2^{l''}\|\dot{\Delta}_{l''} u\|_{L^2}2^{k-l}2^l\|\widetilde{\dot{\Delta}}_lu\|_{L^2}\\
\leq&C\|\rho\|_{L^\infty}\|\dot{S}_{k} u\|_{\dot{B}_{\infty,\infty}^{-1}}\|\nabla u\|_{L^2}\|\nabla\dot{S}_{[\frac{ k}{2}]} u\|_{L^2}.
\end{align*}

For the term $L_3$, we observe that
\begin{align*}
&supp\widehat{\widetilde{\dot{\Delta }}_{l'}\rho\dot{\Delta }_{l'} u_j}\subset\{\xi:|\xi|<5\times2^{l'+1}\},\\
&\dot{\Delta }_l(\widetilde{\dot{\Delta }}_{l'}\rho\dot{\Delta }_{l'} u_j)^k\neq0\Rightarrow l'\geq k-3.
\end{align*}
Similarly, we decompose $L_3$ into two parts.
\begin{align*}
L_3=&\int_{\mathbb{R}^3}\sum\limits_{l\geq k-1}\sum\limits_{l'\geq k-3}\partial_j\dot{S}_k u_i\dot{\Delta }_l\Big(\widetilde{\dot{\Delta }}_{l'}\rho\dot{\Delta }_{l'} u_j\Big)^k\widetilde{\dot{\Delta}}_lu^k_idx\\
=&\int_{\mathbb{R}^3}\sum\limits_{l\geq k-1}\sum\limits_{k-3 \leq l'\leq[\frac{k}{2}]-1 }\partial_j\dot{S}_k u_i\dot{\Delta }_l\Big(\widetilde{\dot{\Delta }}_{l'}\rho\dot{\Delta }_{l'} u_j\Big)^k\widetilde{\dot{\Delta}}_lu^k_idx\\
&+\int_{\mathbb{R}^3}\sum\limits_{l\geq k-1}\sum\limits_{l'\geq [\frac{k}{2}]}\partial_j\dot{S}_k u_i\dot{\Delta }_l\Big(\widetilde{\dot{\Delta }}_{l'}\rho\dot{\Delta }_{l'} u_j\Big)^k\widetilde{\dot{\Delta}}_lu^k_idx\\
:=&L_{31}+L_{32}.
\end{align*}

We estimate the term $L_{32}$ as follows.
\begin{align*}
L_{32}=&\int_{\mathbb{R}^3}\sum\limits_{l\geq k-1}\sum\limits_{l'\geq [\frac{k}{2}]}\partial_j\dot{S}_k u_i\dot{\Delta }_l\Big(\widetilde{\dot{\Delta }}_{l'}\rho\dot{\Delta }_{l'} u_j\Big)^k\widetilde{\dot{\Delta}}_lu^k_idx\\
\leq&\|\nabla\dot{S}_k u\|_{L^\infty}\sum\limits_{l\geq k-1}\sum\limits_{l'\geq [\frac{k}{2}]}\|\widetilde{\dot{\Delta }}_{l'}\rho\|_{L^\infty}\|\dot{\Delta }_{l'} u\|_{L^2}\|\widetilde{\dot{\Delta}}_lu\|_{L^2}\\
\leq& C 2^{\frac{3}{2}k}\|\nabla \dot{S}_k u\|_{L^2}\|\rho\|_{L^\infty}\sum\limits_{l\geq k-1}\sum\limits_{l'\geq [\frac{k}{2}]}2^{l'}\|\dot{\Delta }_{l'} u\|_{L^2}2^l\|\widetilde{\dot{\Delta}}_lu\|_{L^2}2^{-l-l'}\\
\leq& C \|\nabla \dot{S}_k u\|_{L^2}\|\rho\|_{L^\infty}\sum\limits_{l\geq k-1}\sum\limits_{l'\geq [\frac{k}{2}]}2^{(\frac{k}{2}-l')}2^{l'}\|\dot{\Delta }_{l'} u\|_{L^2}2^{(k-l)}2^l\|\widetilde{\dot{\Delta}}_lu\|_{L^2}\\
\leq& C \| \nabla\dot{S}_k u\|_{L^2}\|\rho\|_{L^\infty}\|\nabla u\|_{L^2}^2.
\end{align*}

For the term $L_{31}$, we have
\begin{align*}
L_{31}=&\int_{\mathbb{R}^3}\sum\limits_{l\geq k-1}\sum\limits_{k-4 \leq l'\leq[\frac{k}{2}]-1}\partial_j\dot{S}_k u_i\dot{\Delta }_l\Big(\widetilde{\dot{\Delta }}_{l'}\rho\dot{\Delta }_{l'} u_j\Big)^k\widetilde{\dot{\Delta}}_lu^k_idx\\\
\leq& 2^k\|\dot{S}_k u\|_{L^\infty}\sum\limits_{l\geq k-1}\sum\limits_{k-4 \leq l'\leq[\frac{k}{2}]-1}\|\widetilde{\dot{\Delta }}_{l'}\rho\|_{L^\infty}\|\dot{\Delta }_{l'} u\|_{L^2}\|\widetilde{\dot{\Delta}}_lu\|_{L^2}\\
\leq&C2^{-k}\|\dot{S}_k u\|_{L^\infty}\|\rho\|_{L^\infty}\sum\limits_{l\geq k-1}\sum\limits_{k-4 \leq l'\leq[\frac{k}{2}]-1}2^{(k-l')}2^{l'}\|\widetilde{\dot{\Delta }}_{l'}u\|_{L^2}2^{(k-l)}2^l\|\widetilde{\dot{\Delta}}_lu\|_{L^2}\\
\leq&C2^{-k}\|\dot{S}_k u\|_{L^\infty}\|\rho\|_{L^\infty}\|\nabla u\|_{L^2}\|\nabla\dot{S}_{[\frac{ k}{2}]} u\|_{L^2}\\
\leq&C \|\rho\|_{L^\infty}\|\dot{S}_{k} u\|_{\dot{B}_{\infty,\infty}^{-1}} \|\nabla u\|_{L^2}\|\nabla\dot{S}_{[\frac{ k}{2}]} u\|_{L^2}.
\end{align*}
Combining the conditions \eqref{con1}, $\|\nabla u\|_{L^2}<\infty$ and $\lim\limits_{k\rightarrow-\infty}\|\nabla\dot{S}_{k}u\|_{L^2}=0$,
we conclude that
$$\liminf\limits_{k\rightarrow-\infty}\{ L_1+L_2+L_3\}=0.$$
This completes the proof of Proposition \ref{pro3}.\qquad $\hfill\Box$\\
Consequently, we deduce from Proposition \ref{pro1}-\ref{pro3} and \eqref{t10} that
\begin{align*}
\int_{\mathbb{R}^3}|\nabla u|^2dx=0.
\end{align*}By using the fact
$$\|u\|_{L^6(\mathbb{R}^3)}\leq C\|\nabla u\|_{L^2(\mathbb{R}^3)},$$
we conclude that $u=0$. This completes the proof.\qquad $\hfill\Box$\\

\section*{Acknowledgments}
The authors are supported by the Construct Program of the Key Discipline in Hunan Province and NSFC (Grant No. 11871209), and the Hunan Provincial NSF (No. 2022JJ10033)

\end{document}